\documentclass[12pt]{amsart}
\usepackage{graphicx}

\usepackage[colorlinks=true, pdfstartview=FitV, linkcolor=blue,
citecolor=blue, urlcolor=blue]{hyperref}

\usepackage{mathtools}
\usepackage{amsfonts}
\usepackage{amsmath}
\usepackage{amssymb} 
\usepackage{amsthm}

\usepackage{times}
\usepackage[T1]{fontenc}
\usepackage{enumitem}
\usepackage{setspace}
\usepackage{microtype}
\usepackage{cite}

\allowdisplaybreaks

\usepackage[margin=1.1in]{geometry}

\numberwithin{equation}{section}

\newcommand{\R}{\mathbb R}

\newcommand{\rn}{{{\mathbb R}^n}}

\DeclareMathOperator{\diam}{diam}
\DeclareMathOperator*{\essinf}{ess\,inf}
\DeclareMathOperator*{\esssup}{ess\,sup}

\newcommand{\pp}{{p(\cdot)}}

\newcommand{\Lp}{L^{p(\cdot)}}
\newcommand{\Pp}{\mathcal P}

\newcommand{\qq}{{q(\cdot)}}

\def\XXint#1#2#3{{\setbox0=\hbox{$#1{#2#3}{\int}$}
     \vcenter{\hbox{$#2#3$}}\kern-.5\wd0}}

\newtheorem{thm}{Theorem}[section]

\newtheorem{lem}[thm]{Lemma}

\theoremstyle{definition}

\theoremstyle{remark}
\newtheorem{rem}[thm]{Remark}
\newtheorem{example}[thm]{Example}
\numberwithin{equation}{section}

\begin{document}

\title[On the embedding between   $L^{p(\cdot)}(\Omega)$ and
$L(\log L)^{\alpha}(\Omega)$ space]{On the embedding between  the
  variable Lebesgue space $L^{p(\cdot)}(\Omega)$ and the Orlicz space
  $L(\log L)^{\alpha}(\Omega)$}
\author{David Cruz-Uribe,  OFS,  Amiran Gogatishvili and  Tengiz Kopaliani}

\address{David Cruz-Uribe, OFS\\
	Department of Mathematics, The University of Alabama, Tuscaloosa, AL 35487, USA}
\email{dcruzuribe@ua.edu}

\address{Amiran Gogatishvili \\
Institute of Mathematics of the Czech Academy of Sciences \\
\'Zitna 25 \\
115 67 Prague 1, Czech Republic}
 \email{gogatish@math.cas.cz}

\address{Tengiz Kopaliani \\
Faculty of Exact and Natural Sciences\\
I. Javakhi\-shvili Tbilisi State University\\
 University St. 2\\
 0143 Tbilisi, Georgia}
\email{tengiz.kopaliani@tsu.ge}

\date{June 5, 2024}


\thanks{The first author was partially supported by Simons Foundation
  Travel Support for Mathematicians Grant.  The second author  was  partially supported by the Czech Academy of Sciences RVO: 67985840, by Czech Science Foundation (GA\v CR), grant no: 23-04720S.  The research was in part
  supported by the Shota Rustaveli National Science Foundation   (SRNSF), grant FR 21-12353.}
\begin{abstract}
  We give a sharp sufficient condition on the distribution function,
  $|\{x\in \Omega :\,p(x)\leq 1+\lambda\}|$, $\lambda>0$, of the
  exponent function $p(\cdot): \Omega \to [1,\infty)$ that implies the
  embedding of the variable Lebesgue space $L^{p(\cdot)}(\Omega)$ into
  the Orlicz space $L(\log L)^{\alpha}(\Omega)$,
  $\alpha>0$, where $\Omega$ is an open set with finite Lebesgue
  measure.  As applications of our results, we first give conditions that
  imply the strong differentiation of integrals of functions in
  $L^{p(\cdot)}((0,1)^{n})$, $n>1$.  We then consider the integrability of the
  maximal function on variable Lebesgue spaces, where the exponent
  function $p(\cdot)$ approaches $1$ in value on some part of
  domain.  This  result is an improvement of the result in~\cite{CUF2}.
\end{abstract}
\maketitle

\section{Introduction}

The purpose of this paper is to give a sharp sufficient condition for
the embedding of the variable Lebesgue space
$L^\pp(\Omega)$, where $p_-=1$ and $|\Omega|<\infty$,  into the Orlicz
space $L(\log L)^\alpha(\Omega)$, where $\alpha>0$.   This problem was
originally motivated by attempts to generalize Wiener's integrability
result for the Hardy-Littlewood maximal operator~\cite{Win} to the
scale of variable Lebesgue spaces (see\cite{Has,CUF2,MOS,FM}.

To state our result, we first give the 
definitions and some basic properties of  the spaces we are interested in.   The variable
Lebesgue spaces were introduced by Orlicz in 1931, and have been
studied extensively for the past thirty years.  For more information
and for proofs of the properties we use, see~\cite{MR3026953,MR2790542}.   Let  $\Omega\subset
\mathbb{R}^{n}$ be an open set; we will generally assume that $\Omega$
has finite measure.  Equip it with the  Lebesgue measure.  Let
$\Pp(\Omega)$ consist of all measurable functions
$p(\cdot):\Omega\rightarrow[1,\infty]$.  Given $\pp \in \Pp(\Omega)$
and $E\subset \Omega$, 
let
\[ p_-(E) = \essinf_{x\in E} p(x), \qquad
  p_+(E) = \esssup_{x\in E} p(x).  \]
For brevity, we will write $p_-=p_-(\Omega)$ and $p_+=p_+(\Omega)$.   Note that $p_-\geq 1$; we will frequently assume $p_+<\infty$, that
  is, that the exponent function $\pp$ is bounded.   Define the space $L^{p(\cdot)}(\Omega)$ to be the space of measurable  functions $f$ on $\Omega$  such that for some $\lambda>0,$
$$
\rho_{\lambda}(f)=\int_{\Omega_*}(|f(x)|/\lambda)^{p(\cdot)}\,dx+
\lambda^{-1}\|f\|_{L^\infty(\Omega_\infty)} <\infty,
$$
where $\Omega_\infty =\{ x\in\Omega : p(x) = \infty\}$ and
$\Omega_*=\Omega \setminus \Omega_\infty$.
 The set  $L^{p(\cdot)}(\Omega)$ becomes a Banach space when equpped with
 the Luxemburg norm
 $$
 \|f\|_{p(\cdot)}=\inf\{\lambda>0;\,\,\rho_{\lambda}(f)\leq1\}.
 $$
Define the exponent function $\qq$ pointwise by
 \begin{equation} \label{eqn:dual-exp}
 \frac{1}{p(x)}+\frac{1}{q(x)} = 1, 
\end{equation}
 with the convention that $1/\infty=0$.  Then $L^\qq(\Omega)$ is the
 associate space (also referred to as the K\"othe dual) of
 $L^\pp(\Omega)$.  
 If $p_+<\infty$, then the dual space of $\Lp(\Omega)$ is
 $L^\qq(\Omega)$.

 We now consider Orlicz spaces; for more information, we refer
 to~\cite{MR1113700} and~\cite[Chapter~8]{BeS}.  (For a brief summary,
 see~\cite[Chapter~5]{MR2797562}.)  Let
 $M : [0,\infty)\rightarrow [0,\infty)$ be a Young function, that is, a
 strictly increasing, continuous, convex function on $[0,\infty)$ such
 that $M(0)=0$ and $M(t)/t\rightarrow \infty$ as
 $t\rightarrow \infty$. The Orlicz space $L_{M}(\Omega)$ consists of
 all measurable functions $f$ on $\Omega$ such that the function
 $M(f/\lambda)\in L^{1}(\Omega)$ for some $\lambda>0,$.  It becomes a
 Banach space when  equipped with the norm
$$
\|f\|_{L_{M}}=\inf\left\{\lambda>0;\,\,\int_{\Omega}M(|f(x)|/\lambda)\,dx\leq1\right\}.
$$
Associated to every Young function $M$ is its associate function
$\bar{M}$, defined by
\[ \bar{M}(t) = \sup_{s>0} \{ st- M(s)\}.  \]
This function is also a Young function, and $L_{\bar{M}}$ is the associate space of $L_M$. If $M$ satisfies the
$\Delta_2$ condition--that is, $M(2t)\leq CM(t)$ for all $t\geq
t_0$--then $L_{\bar{M}}$ is the dual space of $L_M$.

If $M(t) =u(t)(\log(e+t))^\alpha$, $\alpha>0$, we denote the Orlicz space
$L_M(\Omega)$  by $L(LogL)^{\alpha}(\Omega)$.  If
$M(t)=\exp(t^\alpha)$, we denote the Orlicz space $L_M$ by
$\exp(L^{\alpha})(\Omega)$.  These spaces are the associate spaces of
one another; moreover,  the dual space of
$L(LogL)^{\alpha}(\Omega)$ is $\exp(L^{1/\alpha})(\Omega)$.  These
spaces are often referred to as Zygmund spaces.  

\medskip

Our main result gives a sharp condition for  the embedding
$L^{p(\cdot)}(\Omega)\subset L(\log L)^{\alpha}(\Omega)$ to hold.

\renewcommand{\labelenumi}{(\alph{enumi})}

\begin{thm}\label{mainembedding}
  Fix an open set $\Omega\subset \rn$, $|\Omega|<\infty$,  and $\pp
  \in \Pp(\Omega)$, $p_+<\infty$.  
  \begin{enumerate}
    
\item Given $\alpha>0$, suppose there exists a  constant $C>1$ such that,
  for all $\lambda>0$ sufficiently small,
\begin{equation} \label{Eq.1.4}
  |\{x\in\Omega:\,\,p(x)\leq1+\lambda\}|
  \leq C^{\frac{1}{\lambda}}\lambda^{\alpha/\lambda}\ln^{-\alpha/\lambda}(1/\lambda).
\end{equation}
Then $L^{p(\cdot)}(\Omega)\subset L(\log L)^{\alpha}(\Omega).$

\item Suppose $\theta : [0,\infty)\rightarrow [0,\infty)$ is an
  increasing, differentiable function such that  $\theta(t)\rightarrow
  \infty$ as $t\rightarrow\infty$, the function
  $\theta(1/\lambda)(1/\lambda)^{-\alpha}(\ln(1/\lambda))^{-\alpha}$
  is increasing for all $\lambda>0$ sufficiently small.  If
$$
\liminf_{\lambda\rightarrow0}
|\{x\in \Omega;:\,p(x)\leq 1+\lambda\}|\theta(1/\lambda)^{-1/\lambda}
\lambda^{-\alpha/\lambda}\ln^{\alpha/\lambda}(1/\lambda)>0,
$$
then $L^{p(\cdot)}(\Omega) \not \subset L(\log L)^{\alpha}(\Omega).$
\end{enumerate}
\end{thm}

As part of  our proof we get an embedding theorem for
the exponential Zygmund spaces.

\begin{thm} \label{thm:exp-zygmund}
Given an open set $\Omega\subset \rn$ and $\pp \in \Pp(\Omega)$.
  Fix $\alpha>0$ and suppose there exists a  constant $C>1$ such that,
  for all $\lambda>0$ sufficiently large,
  \begin{equation} \label{eqn:ez1}
    |\{x\in\Omega : p(x)\geq\lambda\}|
    \leq C^{\lambda}(\lambda)^{-\lambda/\alpha}
    (\ln(\lambda))^{-\lambda/\alpha}
\end{equation}
Then $\exp(L^{\alpha})(\Omega) \subset  L^\pp(\Omega)$.
\end{thm}

Such embeddings were previously considered by  Astashkin and
Masty{\l}o \cite[Theorem~3.1]{AM}, who showed that the embedding
$\exp(L^{2})([0,1])\subset L^{p(\cdot)}([0,1])$ holds if
\eqref{eqn:ez1}
holds with $\alpha=2$.  
They showed that this embedding holds for an exponent function if and
only if the sequence of   Rademacher
functions forms a basic sequence in $L^{p(\cdot)}([0,1])$ that equivalent
to the unit vector basis in $\ell_{2}$.  

\begin{rem}
  The embedding of variable Lebesgue spaces into other spaces that are
  ``infinitesimally close'' to the classical $L^p$ spaces, the
  so-called grand and small Lebesgue spaces, has been investigated by
  the first author, Fiorenza and Guzman~\cite{MR3716507}.
\end{rem}

\medskip

The remainder of this paper is organized as follows.  In
Section~\ref{section:proof} we prove Theorems~\ref{mainembedding}
and~\ref{thm:exp-zygmund}.   Even though the Zygmund spaces are
symmetric spaces (that is, rearrangement invariant spaces) and the
Lebesgue spaces are not, key to our proof is an application of
rearrangement techniques due to~\cite{FR} that let us bridge this difference.
In Section~\ref{section:diff} we give an
application of our results to determine when the basis of rectangles
is a differentiation basis for $\Lp$.  In the classical setting, such
results are closely linked to the boundedness properties of the strong
maximal operator; however, in the variable Lebesgue space setting the
strong maximal operator is never bounded on $\Lp$ unless $\pp$ is
constant~\cite{Kop}.   By using our embedding theorem we can avoid
this problem.  Finally, in Section~\ref{section:max-op} we consider
the problem mentioned above that was the original motivation for the
study of these embeddings:  finding conditions on an exponent function
$\pp$ such that for any ball $B$, $\|Mf\|_{L^1(B)} \leq
\|f\|_{\Lp(\rn)}$.   We will review the known results and show how our
embedding theorem can be used to improve them.

\section{Proof of Theorems \ref{mainembedding} and~\ref{thm:exp-zygmund}}
\label{section:proof}

Hereafter, fix an open set $\Omega$ with $|\Omega|<\infty$.  Given a
measurable function $f:\Omega\rightarrow\mathbb{R}$, define the
distribution function of $f$ to be
$d_{f}(\lambda)=|\{t\in \Omega : |f(t)|>\lambda\}|$, $\lambda\geq 0$
and its decreasing rearrangement by
$f^{\ast}(t)=\inf\{s>0 : d_{f}(s)\leq t\}$, $t>0.$ We say that
functions $f$ and $g$ are equimeasurable if
$f^{\ast}(t)=g^{\ast}(t)$ for all $0<t<|\Omega|,$ or equivalently,
$d_{f}(\lambda)=d_{g}(\lambda)$, $\lambda>0.$

For the proof we need two lemmas.  The first gives an alternative
expression for the norm in $\exp(L^a)(\Omega)$.  For a proof,
see~\cite[Corollary~3.4.28]{EE}.  

\begin{lem} \label{lemma:exp-norm}
For all $\alpha>0$, 
\begin{equation}\label{MarcinkNorm}
  \|f\|_{\exp(L^{\alpha})(\Omega)}
  \approx \sup_{0<t\leq|\Omega|}(\ln(|\Omega|e/t))^{-1/\alpha}f^{\ast}(t).
\end{equation}
\end{lem}

The second is a rearrangment inequality of variable Lebesgue spaces.
For a proof, see~\cite[Theorem~3]{FR}.
 \begin{lem} \label{lem_fr}
 Given $\pp\in \Pp(\Omega)$ and $f\in \Lp(\Omega)$, 
 $$
 \|f\|_{\Lp(\Omega)}\leq (1+|\Omega|)\|f^{\ast}\|_{L^{p^{\ast}(\cdot)}([0,|\Omega|])}.
 $$
 \end{lem}

 \medskip
 
 \subsection*{Proof of Theorem \ref{mainembedding}, Part $(a)$}
 Fix an exponent function $\pp$ that satisfies~\eqref{Eq.1.4}.  Let
 $\qq$ be the dual exponent defined by~\eqref{eqn:dual-exp}.  Then we
 have that 
$$
\{x \in \Omega : p(x)\leq 1+\lambda\}=\{x\in \Omega : q(x)\geq (\lambda+1)/\lambda\};
$$
hence,  the condition~\eqref{Eq.1.4} is, after a change of variables
$\lambda\mapsto 1/\lambda$,  equivalent to
\begin{equation}
    |\{x\in\Omega : q(x)\geq\lambda\}|
    \leq C^{\lambda-1}(\lambda-1)^{-\alpha(\lambda-1)}
    (\ln(\lambda-1))^{-\alpha(\lambda-1)}\label{Eq.3.1}
\end{equation}
 for sufficiently large $\lambda.$
 As we noted above,   $L(\log L)^{\alpha}(\Omega)$ and
 $\exp(L^{1/\alpha})(\Omega)$ are associate spaces.  Therefore, by
 \cite[Proposition~2.10, p.~13]{BeS}, 
$ L^{p(\cdot)}(\Omega)\subset L(\log L)^{\alpha}(\Omega)$ if and only
if $\exp(L^{1/\alpha})(\Omega)\subset L^{q(\cdot)}(\Omega)$.  
Consequently, to complete the proof it is sufficient to prove the
embedding  $\exp(L^{1/\alpha})(\Omega)\subset L^{q(\cdot)}(\Omega)$
assuming that~\eqref{Eq.3.1} holds.  

Without loss of generality, by rescaling we may assume that $|\Omega|=1$.
Fix $x_{0}>1$  sufficiently large and a constant $C_1>1$ such that
$$
C^{x-1}((x-1)\ln(x-1))^{-\alpha(x-1)}\leq C_{1}^{x}(x\ln x)^{-\alpha x}:=F(x)
$$
and the function $F$ is decreasing if $x\geq x_{0}.$ Let $l=F^{-1}$ be
the inverse of $F$ on the interval $(0,t_{0}],$ where
$l(t_{0})=x_{0}.$ From \eqref{Eq.3.1} and the definition of decreasing
rearrangements, we have $q^{\ast}(t)\leq l(t)$ for all
$0<t\leq t_{0}.$ By~\eqref{Eq.1.4}, if
$f\in \exp(L^{1/\alpha})(\Omega)$, then
$f^{\ast}(t)\leq c \ln^{\alpha}(e/t))$.  By Lemma \ref{lem_fr} we need
to prove $\|f^*\|_{L^{p^*(\cdot)}([0,1])}<\infty$; to show this it
  suffices to prove that for some $\lambda>0$,
 $$
\int_0^{t_0} \bigg(\frac{f^*(t)}{\lambda}\bigg)^{p^*(t)}\,dt
\leq
\int_{0}^{t_{0}}\left(\frac{\ln^{\alpha}(e/t)}{\lambda}\right)^{l(t)}\,
dt
=I_\lambda<\infty.
$$
But we can estimate as follows:  making the change of variables $x
\mapsto F(x)$,
\begin{align*}
  I_{\lambda}
  & =-\int_{x_{0}}^{\infty}(\lambda^{-1}(\ln(eC_{1}^{-x}(x\ln
    x)^{\alpha x}))^{\alpha})^{x}
    \,dF(x) \\
  & =\alpha\int_{x_{0}}^{\infty}\lambda^{-x}(\alpha x)^{\alpha x}
    \ln^{\alpha x}(e^{\frac{1}{\alpha x}}C_{1}^{-1/\alpha}x\ln x) \\
  & \qquad \times C_{1}^{x}(x\ln x)^{-\alpha x}
    \left(\ln(C_{1}^{-\frac{1}{\alpha}}x\ln x)
    +\frac{\ln x+1}{\ln x}\right)\,dx \\
& \leq C_{2}\int_{x_{0}}^{\infty}\left(\frac{C_{3}}{\lambda}\right)^{x}
 \left(\ln(x\ln x)+\frac{\ln x+1}{\ln x}\right)\,dx \\
&   <\infty;
\end{align*}
the last inequality holds provided that we choose $\lambda>C_{3},$.
This completes the proof.

\medskip

\subsection*{Proof of Theorem~\ref{thm:exp-zygmund}}
Our hypothesis~\eqref{eqn:ez1} can be rewritten as
\[ |\{ x\in \Omega : p(x) \geq \lambda \} \leq G(x), \]
where $G(x) = C^x(x\ln(x)|^{-x/\alpha}$.  Therefore, we can repeat the
above argument, replacing $F$ by $G$ (and so $\alpha$ by $1/\alpha$), and we get the desired embedding.

\medskip

\subsection*{Proof of Theorem \ref{mainembedding}, Part $(b)$}
By our hypothesis, there exists $\lambda_0,\,\gamma\in (0,1)$ such that for all $0<\lambda\leq \lambda_{0}$ we have
$$
|\{x\in \Omega : p(x)\leq1+\lambda\}|\geq\gamma(\theta(1/\lambda)\lambda{\alpha}(\ln(1/\lambda))^{-\alpha})^{1/\lambda}.
$$
Arguing as before, this inequality is equivalent to following
inequality for the dual exponent $\qq$:
\begin{equation}
|\{x \in \Omega :\,q(x)\geq\lambda\}|\geq\gamma(\theta(\lambda-1)((\lambda-1)\ln(\lambda-1))^{-\alpha})^{\lambda-1},\label{Eq.3.2a}
\end{equation}
for all $\lambda\geq1+1/\lambda_{0}$.  Again by duality, it will
suffice to show that
$\exp(L^{1/\alpha})(\Omega) \not\subset L^{q(\cdot)}(\Omega).$

Without loss of generality we will assume that $|\Omega|=1$.  Note
that   given a measure preserving transformation
$\omega:\,\Omega\rightarrow(0,1]$, for all $t\in (0,1]$,
$$
\left(\ln^{\alpha}\left(\frac{e}{\omega(\cdot)}\right)\right)^{\ast}(t)
=\ln^{\alpha}(e/t);
$$
hence, by Lemma~\eqref{MarcinkNorm}, 
$\ln^{\alpha}(e/\omega(x))\in \exp(L^{1/\alpha})(\Omega)$.
We will show that for every $c>0$ there exists a measure preserving transformation $\omega:\,\Omega\rightarrow(0,1]$
such that
\begin{equation}
\int_{\Omega}(c^{-1}\ln^{\alpha}(e/\omega(x)))^{q(x)}dx=\infty. \label{Eq.3.2}
\end{equation}
It follows at once from (\ref{Eq.3.2}) and the definition of the norm
that  $\|\ln^{\alpha}(e/\omega)\|_{q(\cdot)}>c$; consequently, the
inequality  $\|f\|_{q(\cdot)}\leq C\|f\|_{\exp(L^{1/\alpha})}$  cannot
hold for every $f\in \exp(L^{1/\alpha})$  with $C$ is independent of
$f$.

To prove \eqref{Eq.3.2}, fix $c>0$ large.   Without loss of
generality, we may assume that for all $\lambda$ sufficiently large,
\begin{equation}  \label{Eq.3.3}
\theta(\lambda-1)\leq \ln^{\alpha}(\lambda-1);
\end{equation}
otherwise, instead of $\theta(\lambda-1)$ we can take the function $\min\{\theta(\lambda-1),\ln^{\alpha}(\lambda-1)\}.$ Moreover, our hypotheses on $\theta$ imply
\begin{align*}
\left(\frac{\theta^{\frac{1}{\alpha}}(\lambda-1)}{(\lambda-1)\ln
  (\lambda-1)}\right)'
  & =
\theta^{\frac{1}{\alpha}}(\lambda-1)(\lambda-1)^{-2}\ln^{-2}(\lambda-1) \\
  & \qquad \times \left(\frac{\theta'(\lambda-1)}
    {\alpha\theta(\lambda-1)}(\lambda-1)\ln (\lambda-1)
    -(1+\ln (\lambda-1))\right) \\
  & \leq 0,
\end{align*}
and so for $\lambda\geq \widetilde{\lambda}_{0}$,
\begin{equation}
\frac{(\lambda-1)\theta'(\lambda-1)}{\alpha\theta(\lambda-1)}\leq\frac{1+\ln (\lambda-1)}{\ln(\lambda-1)},\,\,\,, \label{Eq.3.4}
\end{equation}
provided we fix $\widetilde{\lambda}_{0}$ sufficiently large.

Define $\psi(\lambda)=\theta(\lambda)(\lambda\ln(\lambda))^{-\alpha}$,
and let  $g$ be the inverse function to the function $\gamma
\psi(\lambda-1)^{\lambda-1}$ when $\lambda\geq \widetilde{\lambda}_{0}$.
It follows from \eqref{Eq.3.2a} that for some $t_{0}\in(0,1].$
\begin{equation}
q^{\ast}(t)\geq g(t),\,\,\,0<t\leq t_{0},  \label{Eq.3.5}
\end{equation}
 We may  assume that, for a given $c,$
the inequality $\ln^{\alpha}(e/t)\geq c$ is valid for all $t\in(0,t_{0}].$

By Ryff's theorem \cite[Theorem~7.5, p. 82]{BeS},  there exists a
measure preserving transformation $\omega:\,\Omega\rightarrow(0,1]$
such that $q(x)=q^{\ast}(\omega(x))$ almost everywhere.   If we define
$E=\omega^{-1}([0,t_{0}])$, then it follows from  inequality
(\ref{Eq.3.5})  that for almost every $x\in E$,
$$
q(x)\geq g(\omega(x)).
$$
Therefore, 
\begin{multline*}
  I_{c}=\int_{E}(c^{-1}\ln^{\alpha}(e/\omega(x)))^{q(x)}\,dx \\
  \geq\int_{E}(c^{-1}\ln^{\alpha}(e/\omega(x)))^{g(\omega(x))}\,dx
=\int_{0}^{t_{0}}(c^{-1}\ln^{\alpha}(e/t))^{g(t)}\,dt;
\end{multline*}
if we make the change of variables $\lambda=g(t)$, we obtain
$$
I_{c}\geq -\gamma\int_{g(t_{0})}^{\infty}c^{-\lambda}\ln^{\alpha \lambda}\left(\frac{e}{\gamma(\psi(\lambda-1))^{\lambda-1}}\right)d((\psi(\lambda-1))^{\lambda-1}).
$$

To estimate this integral, note first that 
\begin{multline*}
  ((\psi(\lambda-1))^{\lambda-1})'
  =\left(\exp\left(-\alpha
      (\lambda-1)\ln(\theta^{-\frac{1}{\alpha}}(\lambda-1)(\lambda-1)\ln
      (\lambda-1))\right)\right)' \\
=-\alpha(\psi(\lambda-1))^{\lambda-1}\left(\ln\frac{(\lambda-1)\ln (\lambda-1)}{\theta^{1/\alpha}(\lambda-1)}+\frac{1+\ln (\lambda-1)}{\ln (\lambda-1)}-\frac{(\lambda-1)\theta'((\lambda-1))}{\alpha\theta(\lambda-1)}\right).
\end{multline*}
But then, by (\ref{Eq.3.4}), we have that for
$\widetilde{\lambda}_{0}$ sufficiently large and all
$\lambda\geq \widetilde{\lambda}_{0}$,
$$
((\psi(\lambda-1))^{\lambda-1})'
\leq-\alpha\psi(\lambda-1)^{\lambda-1}\ln\frac{(\lambda-1)
  \ln(\lambda-1)}{\theta^{1/\alpha}(\lambda-1)}.
$$
Since
$(\theta^{-1/\alpha}(1/\lambda)(1/\lambda)(\ln(1/\lambda)))^{-\alpha}$
is increasing for small enough $\lambda$,  we have that 
$\theta^{-\frac{1}{\alpha}}(\lambda-1)(\lambda-1)\ln (\lambda-1))$ is
increasing for large $\lambda$.  If we combine this with inequality
(\ref{Eq.3.3}), we obtain (for is sufficiently large
$\widetilde{\lambda_{0}}$)
\begin{align*}
  I_{c}
  & \geq \gamma\alpha\int_{g(t_{0})}^{\infty}c^{-\lambda}
    (\alpha (\lambda-1))^{\alpha \lambda}\ln^{\alpha \lambda}
    (\gamma^{-\frac{1}{\alpha(\lambda-1)}}e^{\frac{1}{\alpha(
    \lambda-1)}}
    \theta^{-\frac{1}{\alpha}}(\lambda-1)(\lambda-1)\ln (\lambda-1)) \\
& \qquad \times \theta^{\lambda-1}(\lambda-1)(\lambda-1)^{-\alpha
 (\lambda-1)} \ln^{-\alpha (\lambda-1)}(\lambda-1)
 \ln\frac{(\lambda-1)\ln (\lambda-1)}
 {\theta^{\frac{1}{\alpha}}(\lambda-1)}\,d\lambda \\
  & \geq \gamma\alpha\int_{g(t_{0})}^{\infty}
    (c\alpha^{-\alpha})^{-\lambda}
    \theta^{\lambda-1}(\lambda-1)\ln(\lambda-1)\,d\lambda.
\end{align*}

Since $\theta(x) \rightarrow \infty$ as $x\rightarrow\infty$,   it
dominates the exponential term in the last integral, and it follows
that $I_{c}=\infty$, which implies (\ref{Eq.3.2}).  This completes the proof.

\section{Strong differentiation of functions from  variable Lebesgue
  spaces $L^{p(\cdot)}((0,1)^n)$}
\label{section:diff}

In this section we apply Theorem~\ref{mainembedding} to give a
sufficient condition on the exponent $\pp$ for the basis of rectangles
to be a differentiation basis of $\Lp(\Omega)$.  We first recall some
basic definitions and results.  The family
$\mathcal{B}=\{B(x):x\in\Omega\}$ is called a basis in $\Omega$ if,  for
each $x\in \Omega$,  $B(x)$ is a collection of bounded, open sets that
contain
$x$ and there exists a sequence $\{B_{k}\}\subset B(x),$ such that
$\diam(B_{k})\rightarrow0$ as $k\rightarrow\infty.$

Given $f\in L^1_{loc},$ the numbers
\begin{align*}
  \overline{D}_{\mathcal{B}}\left(\int f,x\right)
  & =\sup_{\{\{B_{k}\}\subset
  B(x)\}}
\limsup_{\diam(B_{k})\rightarrow0}|B_{k}|^{-1}\int_{B_{k}}f(y)\,dy, \\
  \underline{D}_{\mathcal{B}}\left(\int f,x\right)
  & =\inf_{\{\{B_{k}\}\subset
    B(x)\}}\liminf_{\diam (B_{k})\rightarrow0}
    |B_{k}|^{-1}\int_{B_{k}}f(y)\,dy
\end{align*}
    are called, respectively, the upper and the lower derivatives of
    the integral of $f$ at the point $x.$ If the upper and lower
    derivatives coincide, then their common value is called the
    derivative of the integral of $f$ at the point $x$ and denoted
    by $D_{\mathcal{B}}\left(\int f,x\right).$ A basis $\mathcal{B}$
    is said to differentiate the integral of $f$ if
    $D_{\mathcal{B}}\left(\int f,x\right)=f(x)$ for almost all $x.$

    We are interested in the basis of rectangles, that is, the basis
    $\mathcal{B}_{1}$ where for each $x$, $B(x)$ contains all
    rectangles with the sides parallel to the axes containing $x$ and
    contained in $\Omega$. For simplicity, let
    $\Omega=\mathbb{T}^{n}=(0,1)^{n},(n\geq1).$ By the
    Jessen-Marcinkiewicz-Zygmund theorem \cite{JMZ} (see
    also~\cite[p.~50]{MR0457661}), for $n>1$ the basis
    $\mathcal{B}_{1}$ differentiates functions from
    $L(\log L)^{n-1}(\mathbb{T}^{n})$.  On the other hand Saks~\cite{Sak} proved
    that for each function $\sigma(t)$ that decreases monotonically to
    zero, there exists a function
    $f\in \sigma(L)L(\log L)^{n-1}(\mathbb{T}^{n}),(n>1)$ such that
    $\overline{D}_{\mathcal{B}_{1}}\left(\int f,x\right)=+\infty$
    a.e. on $\mathbb{T}^{n}$.

    Closely related to the differentiation properties of
    $\mathcal{B}_{1}$ are the boundedness properties of  $M_{S}f$, the
    maximal operator  corresponding to $\mathcal{B}_{1}$, also known
    as the strong maximal function:
$$
M_{S}f(x)=\sup_{x\in I\in \mathcal{B}_{1}}|I|^{-1}\int_{I}|f(y)|dy.
$$
The strong maximal operator is bounded on $L^{p}(\mathbb{T}^{n}),$ for
$p>1$~\cite{JMZ}, but is not bounded on
$L^{p(\cdot)}(\mathbb{T}^{n})$, $n>1$, if $p(\cdot)$ is not
constant~\cite{Kop}.  Now if $p_->1$, then
$ L^{p(\cdot)}(\mathbb{T}^{n})\subset L^{p_{-}}(\mathbb{T}^{n})$ and
consequently, if  $f\in L^{p(\cdot)}(\mathbb{T}^{n})$,
$D_{\mathcal{B}_{1}}\left(\int f,x\right)=f(x)$ a.e.  We are
interested in finding conditions on $\pp$ for which 
$p(x)>1$ a.e. on $\mathbb{T}^{n}$ but $p_-=1$, and indeed $p_-(E)=1$
for any open set $E\subset \mathbb{T}^{n}$.  
For this problem it is natural find  conditions on the function
$|\{x:p(x)\leq 1+\lambda\}|$, $\lambda>0$, that  guarantee the differentiation of functions from $ L^{p(\cdot)}(\mathbb{T}^{n}),\,(n>1).$ 
From part $(a)$ of Theorem~\ref{mainembedding} and the
Jessen-Marcinkiewicz-Zigmund theorem   we immediately obtain following
result.

\begin{thm}\label{main}
 Given $\pp \in \Pp(\mathbb{T}^{n})$, $n>1$, suppose that  for some $C>1$ and
 for all $\lambda>0$ sufficiently small, 
 $$
 |\{x\in\mathbb{T}^{n}:p(x)\leq1+\lambda\}|\leq C^{1/\lambda}\lambda^{(n-1)/\lambda}(\ln(1/\lambda))^{-(n-1)/\lambda}.
 $$
 Then for $f\in L^{p(\cdot)}(\mathbb{T}^{n})$ and almost every
 $x\in \mathbb{T}^{n}$,
 $D_{\mathcal{B}_{1}}\left(\int f,x\right)=f(x)$.
\end{thm}

\begin{rem}
The case when exponent $p(x)>1$ a.e, and the set $\{x:p(x)=1\}$ is
"irregular" has other interesting applications. For instance,  in
\cite{EGK} the authors constructed  an exponent $\pp \in
\Pp([0,2\pi])$ such that  $p(x)>1$ a.e., and 
the space $L^{p(\cdot)}[0,2\pi]$  contains both the Kolmogorov and the
Marcinkiewicz examples of functions in $L^{1}$ with a.e. divergent
Fourier series.
\end{rem}

\medskip

We can also give examples of exponents for which $\mathcal{B}_1$ fails
to be a differentiation basis for $\Lp(\mathbb{T}^{n})$.  
Stokolos \cite{Stok} in case $n=2$  and Oniani \cite{ON} in case $n>2$
proved following result.

\begin{thm}\label{StokolosOniani}
For each function $f\in L^1(\mathbb{T}^{n})$ such that $f\not \in L(\log L)^{n-1}(\mathbb{T}^{n}),$ there exists a function $g$ that is equimeasurable with $f$ such that $\overline{D}_{\mathcal{B}_{1}}\left(\int f,x\right)$ $=+\infty$ a.e.
\end{thm}

In fact the function $g$ in Theorem~\ref{StokolosOniani} has the form
$g(x)=f(\omega(x))$, where $\omega:\mathbb{T}^{n}\rightarrow
\mathbb{T}^{n}$ is an invertible, measure-preserving
transformation. Therefore, from  part $(b)$ of
Theorem~\ref{mainembedding} and Theorem~\ref{StokolosOniani} we obtain
get the following result.

\begin{thm}
For $n>1$, let $\pp \in \Pp(\mathbb{T}^{n})$ be such that  for some $0\leq \epsilon<n-1$,
$$
|\{x\in\mathbb{T}^{n}:p(x)\leq1+\lambda\}|\geq C \lambda^{\epsilon/\lambda}(\ln(1/\lambda))^{-(n-1)/\lambda};
$$
alternatively, let $\pp$ be such htat 
$$
|\{x\in\mathbb{T}^{n}:p(x)\leq1+\lambda\}|\geq C
\lambda^{(n-1)/\lambda}(\ln(1/\lambda))^{-\epsilon/\lambda}.
$$
Then there exists $\bar{p}(\cdot) \in \Pp(\mathbb{T}^{n})$
equimeasurable with $\pp$,  and   $f\in L^{\bar{p}(\cdot)}(\mathbb{T}^{n})$  such that $\overline{D}_{\mathcal{B}_{1}}\left(\int f,x\right)=+\infty$ a.e.
\end{thm}

For specific examples of such exponent functions, see Examples~\ref{ex:one}
and~\ref{ex:three} below.

\section{Further remarks on the embedding $L^{p(\cdot)}(\Omega)\subset
  L(\log L)^{\alpha}(\Omega)$.}
\label{section:max-op}

In this section we consider one of the original motivations for
studying the embedding $L^{p(\cdot)}(\Omega)\subset
  L(\log L)^{\alpha}(\Omega)$:  the generalization of Weiner's theorem
  to variable Lebesgue spaces.  We recall some background
  information.  Given an open set $\Omega\subset\mathbb{R}^{n}$ and
  $f\in L^1_{loc}(\Omega)$, we define the Hardy-Littlewood maximal function by
$$
Mf(x)=\sup_{x\in B}\frac{1}{|B|}\int_{B\cap\Omega}|f(y)|\,dy,
$$
where the supremum is taken over balls containing $x.$
It is well known that for $1<p<\infty$, $M:\,L^{p}(\Omega)\rightarrow
L^{p}(\Omega),$; on the other hand, given any $f\in
L^{1}(\mathbb{R}^{n})$, $f\neq 0,$  $Mf \not\in L^{1}(\mathbb{R}^{n})$. In fact, $Mf$ need not even be locally integrable.
Wiener \cite{Win} proved that $Mf$ is locally integrable if $f \in L\log L(\mathbb{R}^{n})$. More precisely, he
showed that given any ball $B$,
\begin{equation}\label{Winer}
\|Mf\|_{L^{1}(B)}\leq C\|f\|_{L\log L(\mathbb{R}^{n})}.
\end{equation}

The question arose to determine conditions on an exponent $\pp\in
\Pp(\rn)$ so that for a set $\Omega \in \rn$ with finite measure,
\begin{equation} \label{eqn:var-wiener}
\|Mf\|_{L^{1}(\Omega)}\leq C\|f\|_{p(\cdot)}. 
\end{equation}
In particular, if $p_-=1$, determine when~\eqref{eqn:var-wiener}
holds.

The first result in this direction was due to  H\"{a}st\"{o}~\cite{Has}.
Fix numbers $a>0$ and $b\in \R$, and define the function
\begin{equation} \label{Eq.1.1}
\Lambda(r)=\Lambda(r,a,b)=1+\frac{a\ln\ln(1/r)}{\ln(1/r)}+\frac{b}{\ln(1/r)} 
\end{equation}
for $0<r<r_{0}$ and $\Lambda(r)=\Lambda(r_{0})$ for $r>r_{0}.$
For a compact set $K$ in $\mathbb{R}^{n},$  define
$$
K(r)=\{x\in\mathbb{R}^{n}:\,\,\delta_{K}(x)< r\},
$$
where $\delta_{K}(x)$
denotes the distance
from $x$ to set $K.$ For $\beta>0,$ we say that the Minkowski
$(n-\beta)$ content of $K$ is finite if $|K(r)|\leq Cr^{\beta}$ for
all $r>0$ sufficiently small.
Now we define
a variable exponent $p(\cdot)$ by
\begin{equation}
p(x)=\Lambda(\delta_{K}(x)) \label{Eq.1.2}
\end{equation}
for $x\in \mathbb{R}^{n};$ set $p(x)=1$ on $K.$

\begin{thm}\label{Has}
 Let $K$ be a compact subset of a bounded, open set $\Omega$,  and let
 $p(\cdot)$ be a given by (\ref{Eq.1.2}). If $a>1$ and the Minkowski
 $(n-1)$-content of $K$ is finite, then~\eqref{eqn:var-wiener} holds.
\end{thm}

The heart of the proof of Theorem~\ref{Has} is to show that given
these hypotheses, $L^{p(\cdot)}(\Omega)\subset L\log L(\Omega)$.  The
desired conclusion then follows immediately from
inequality~\eqref{Winer}.  Later, Futamura and Mizuta \cite{FM} proved
that the conclusion is still valid when $a=1$ and $b=0$.

Mizuta, Ohno, and Shimomura~\cite{MOS} proved the following embedding theorem.

\begin{thm}\label{MOS} Let $K$ be a compact subset of bounded open set
  $\Omega$ whose Minkowski $(n-\beta)$-content is finite. Set
  $p(x)=\Lambda(\delta_{K}(x))$, $x\in \Omega$. If $a>0$, then
  $L^{p(\cdot)}(\Omega)\subset L(\log L)^{a\beta}(\Omega).$
\end{thm}

If we set $a=1/\beta$  in Theorem \ref{MOS}, we get a generalization
of Theorem \ref{Has}.  

Subsequently, the first author and Fiorenza \cite{CUF2} proved the
following result.

\begin{thm}\label{GU} Given $ \pp \in \Pp(\rn)$,   suppose there exist
  constants $\epsilon$, $0<\epsilon<1$, $C_0>0$,  and
  $\delta$, $0<\delta<e^{-e},$ such that for $0<\lambda\leq\delta,$
\begin{equation}
  |\{x\in \rn : p(x)\leq \Lambda(\lambda,1+\epsilon,0)\}|
  \leq C_0\lambda.\label{Eq.1.3}
\end{equation}
Then given any ball $B$, there exists a constant $C$ (depending on
$|B|$, $p(\cdot)$, $\epsilon$, and $\delta$) such
that~\eqref{eqn:var-wiener} holds (with $\Omega$ replaced by $B$).
\end{thm}

One feature of their proof is that they did not pass through the
embedding $L^{p(\cdot)}(\Omega)\subset L\log L(\Omega)$; rather they
gave a direct proof of~\eqref{eqn:var-wiener}.  Here we will show that
we can use Theorem~\ref{mainembedding} to give a small improvement of
Theorem~\ref{GU} which is analogous to the results of Futamura and
Mizuta, and Mizuta, {\em et al.}

To do so, we will use the Lambert $W$-function to find a relationship
between estimates of the form (\ref{Eq.1.3}) and (\ref{Eq.1.4}).
The Lambert $W$-function is the solution of the equation
$$
W(x)\exp(W(x))=x.
$$
Here we summarize some of its properties; for complete information see
\cite{BeoC,CGHJK}.  Note that $W(-1/e)=-1$.  On the interval
$[0,\infty)$ there is one real solution and on the interval $(-1/e,0)$
there are two real solutions. Call the solution for which
$W(x)\geq W(-1/e)$ the principal branch and denote it by $W_p(x).$
This function is defined on $[-1/e;\infty)$ and
$W_{p}:[-1/e;+\infty)\rightarrow[-1,+\infty)$ is increasing bijection.
Denote the other solution  $W_m(x).$ This function is defined on
$(-1/e;0)$ and $W_{m}:(-1/e;0)\rightarrow[-\infty,-1)$ is a decreasing
bijection.  Moreover, we have that as $x\rightarrow+\infty$,
\begin{equation}
W_p(x)=\xi-\ln\xi+\frac{\ln\xi}{\xi}+\frac{(\ln\xi)^{2}}{2\xi^{2}}-\frac{(\ln\xi)}{\xi^{2}}+O\left(\frac{(\ln\xi)^{3}}{\xi^{3}}\right),\label{Eq.1.51}
\end{equation}
where $\xi=\ln x$.   Similarly, we have that as $x\rightarrow 0$, 
\begin{equation}
W_m(x)=-\mu-\ln\mu-\frac{\ln\mu}{\eta}-\frac{(\ln\mu)^{2}}{2\mu^{2}}-\frac{(\ln\mu)}{2\mu^{2}}+O\left(\frac{(\ln\mu)^{3}}{\mu^{3}}\right),\label{Eq.1.5}
\end{equation}
where $\mu=\ln(-1/x)$.

Fix $a,b\in \R$ with $a>0.$ For $0<r\leq r_{0}$ ($r_{0}$ a small
positive number) we will find the inverse function of the function
$\Lambda(r,a,b)-1$.   The equation $\Lambda(r,a,b)-1=x$ may
rewritten as
$$
-\ln(k\ln(1/r))\exp(-\ln(k\ln(1/r)))=-\frac{x}{ak},
$$
where $k=\exp(b/a).$
Since $-\ln(k\ln(1/r))\rightarrow-\infty$ as $r\rightarrow 0^+$, for
$r>0$ small we have that
$$
r=\exp\left(-\frac{1}{k}\exp\left(-W_m\left(-\frac{x}{ak}\right)\right)\right).
$$
Using the fact that $\exp(W(x))=x/W(x)$, we obtain
\begin{equation}\label{inverse}
r=\exp\left(\frac{a}{x}W_m\left(-\frac{x}{ak}\right)\right)=
\left(\frac{x}{ak}\right)^{\frac{a}{x}}\left(-W_m\left(-\frac{x}{ak}\right)\right)^{-\frac{a}{x}}.
\end{equation}

If we apply Theorem \ref{mainembedding} using \eqref{inverse},
\eqref{Eq.1.5},  we get the following result.

\begin{thm}
 Fix  $\Omega \subset \rn$, $|\Omega|<\infty$, and $\pp \in
 \Pp(\Omega)$, $p_+<\infty$.   Given $a>0$,   suppose there exists a
 constant $C>1$ such that, 
 for all $\lambda>0$ sufficiently small,
\begin{equation}
|\{x\in\Omega:\,\,p(x)\leq1+\Lambda(\lambda,a,b)\}|\leq C\lambda.
\end{equation}
Then $L^{p(\cdot)}(\Omega)\subset L(\log L)^{a}(\Omega).$  In
particular, if $a=1$ and $b\in \R$, then
inequality~\eqref{eqn:var-wiener} holds.
\end{thm}

We can see directly that the hypothesis of Theorem~\ref{GU} is
stronger than that of Theorem~\ref{mainembedding}.  If we use
\eqref{inverse} with $a=1+\epsilon$ and $b=0$, and then apply
estimate(\ref{Eq.1.5}), we can rewrite (\ref{Eq.1.3}) as
\begin{multline*}
|\{x\in\Omega : p(x)\leq 1+\lambda\}| \\
\leq C_0
\left(\frac{\lambda}{1+\epsilon}\right)^{\frac{1+\epsilon}{\lambda}}
\left(-W_m\left(-\frac{\lambda}{1+\epsilon}\right)
\right)^{-\frac{1+\epsilon}{\lambda}} 
\leq
C\lambda^{\frac{1+\epsilon}{\lambda}}(\ln(1/\lambda))^{-\frac{1+\epsilon}{\lambda}}.
\end{multline*}
In other words,   if for some $\epsilon>0$, \eqref{Eq.1.3} holds,  then
$$
|\{x\in\Omega:\,p(x)\leq 1+\lambda\}|\leq C\lambda^{\frac{1+\epsilon}{\lambda}}(\ln(1/\lambda))^{-\frac{1+\epsilon}{\lambda}}
$$
for sufficiently small $\lambda.$
However, we also have that
$$
\frac{C^{\frac{1}{\lambda}} \lambda^{\frac{1}{\lambda}}(\ln(1/\lambda))^{-\frac{1}{\lambda}}}{\lambda^{\frac{1+\epsilon}{\lambda}}(\ln(1/\lambda))^{-\frac{1+\epsilon}{\lambda}}}\rightarrow+\infty\,\,\,\mbox{as}\,\,\lambda\rightarrow0.
$$
Thus, the hypotheses of Theorem~\ref{mainembedding} are weaker than
those of Theorem~\ref{GU}.

\medskip

We conclude this section with three examples of exponents for which
parts $(a)$ or $(b)$ of Theorem~\ref{mainembedding} hold for $\Omega
\subset \R$, $\Omega=[0,x_0]$ for some $x_0>0$.

\begin{example} \label{ex:one}
  Let $\alpha>0,0\leq\epsilon<\alpha$, $\theta(x)=(\ln
  x)^{\alpha-\epsilon},$ and suppose the exponent  $p(\cdot)$ satisfies
 $$
 |\{x: p(x)\leq 1 +\lambda\}|\geq C\lambda^{\alpha/\lambda}\ln^{-\epsilon/\lambda}(1/\lambda).
 $$
 Then the conditions part $(b)$ of Theorem \ref{mainembedding}  are
 valid,  and so $L^{p(\cdot)}(\Omega)\not \subset L(\log L)^{\alpha}(\Omega).$
\end{example}

We further consider the case $\epsilon=0$ in Example~\ref{ex:one}.  By
using the Lambert $W$-function we may find (arguing as we did above)
the inverse function of $\lambda^{\alpha/\lambda}$ on
$0<\lambda\leq \lambda_{0},$ for sufficiently small $\lambda_{0}.$ We
have that if
$\lambda^{\alpha/\lambda}=x,\,0<x\leq\lambda_{0}^{\alpha/\lambda_{0}}=x_{0},$
then
$$
\ln(1/\lambda)\exp(\ln(1/\lambda))=\frac{1}{\alpha}\ln(1/x).
$$
Again using the fact that $\exp(W(x))=x/W(x)$, we obtain
$$
\lambda=\frac{\alpha W_p\left(\frac{1}{\alpha}\ln(1/x)\right)}{\ln(1/x)}:= p_{0}(x).
$$
Consequently,  for exponent $p(x)=1+p_{0}(x)$ we have that for
$0<\lambda\leq \lambda_{0}$, 
 $$
 |\{x\in[0,x_{0}] : p(x)\leq 1+\lambda\}|=\lambda^{\alpha/\lambda}, 
 $$
 and  $L^{p(\cdot)}([0,x_{0}])\not \subset L(\log L)^{\alpha}([0,x_{0}]).$
 If we apply  (\ref{Eq.1.51}),  we obtain  following expression for $p(x)$:
$$
p(x)=1+\frac{\alpha\ln\left(\frac{1}{\alpha}\ln\left(\frac{1}{x}\right)\right)}{\ln\left(\frac{1}{x}\right)}-
\frac{\alpha\ln\ln\left(\frac{1}{\alpha}\ln\left(\frac{1}{x}\right)\right)}{\ln\left(\frac{1}{x}\right)}+
o\left(\frac{\alpha\ln\ln\left(\frac{1}{\alpha}\ln\left(\frac{1}{x}\right)\right)}{\ln\left(\frac{1}{x}\right)}\right)
$$
as $x\rightarrow 0.$

\begin{example}
For sufficiently small $x_0$, define for $0<x<x_0$ the exponent function
$$
p(x)=1+\frac{\alpha\ln\left(\frac{1}{\alpha}\ln\left(\frac{1}{x}\right)\right)}{\ln\left(\frac{1}{x}\right)}
=1+\frac{\ln\left(\ln\left(\frac{1}{x^{1/\alpha}}\right)\right)}{\ln\left(\frac{1}{x^{1/\alpha}}\right)},\,\,(\alpha>0).
$$
Then, using \eqref{inverse}, we have that for $0<\lambda\leq \lambda_{0}$,
 $$
|\{x\in[0,x_{0}]:\,p(x)\leq 1+\lambda\}|=\lambda^{\alpha/\lambda}(-W_{m}(-\lambda))^{-\alpha/\lambda}\leq C
\lambda^{\alpha/\lambda}\ln^{-\alpha/\lambda}(1/\lambda);
$$
hence,
$L^{p(\cdot)}([0,x_{0}])\subset L(\log L)^{\alpha}([0,x_{0}]).$
\end{example}

\begin{rem}
 Let $K$ be a compact subset of $\Omega$ such that
$$
C^{-1}r^\alpha\leq |K(r)|\leq Cr^{\alpha}\,\,\,\mbox{for}\,\,0<r<r_{0}\,\,\,(\alpha>0)
$$
and define
$$
p(x)=1+\frac{b\log\log(1/\delta(x))}{\log(1/\delta(x))}-\frac{c\log\log\log(1/\delta(x))}{\log(1/\delta(x))},
$$
where $b>0,c>0,\,\delta(x)\leq r_{0}$ and $\inf_{\{x,\,\delta(x)>r_{0}\}}p(x)>1.$
In \cite{MOS} (see also \cite{FM}) the authors  constructed a function
$f\in L^{p(\cdot)}(\Omega)$ such that $f \not \in L(\log L)^{\alpha b}(\Omega).$
\end{rem}

\medskip

\begin{example}\label{ex:three}
  Let $\alpha>0,0\leq\epsilon<\alpha$, $\theta(x)=
  x^{\alpha-\epsilon},$ and suppose exponent function $p(\cdot)$ satisfies
 $$
 |\{x: p(x)\leq 1 +\lambda\}|\geq C\lambda^{\epsilon/\lambda}\ln^{-\alpha/\lambda}(1/\lambda).
 $$
 Then the conditions of part $(b)$ of Theorem \ref{mainembedding} are
 valid, and so  $L^{p(\cdot)}(\Omega) \not \subset L(\log L)^{\alpha}(\Omega).$
\end{example}

\bibliographystyle{plain}
\bibliography{STRONG-DIFF}

\end{document}